# ADAPTIVE VARIANCE FUNCTION ESTIMATION IN HETEROSCEDASTIC NONPARAMETRIC REGRESSION

By T. Tony Cai[1] and Lie Wang

*University of Pennsylvania and University of Pennsylvania*

We consider a wavelet thresholding approach to adaptive variance function estimation in heteroscedastic nonparametric regression. A data-driven estimator is constructed by applying wavelet thresholding to the squared first-order differences of the observations. We show that the variance function estimator is nearly optimally adaptive to the smoothness of both the mean and variance functions. The estimator is shown to achieve the optimal adaptive rate of convergence under the pointwise squared error simultaneously over a range of smoothness classes. The estimator is also adaptively within a logarithmic factor of the minimax risk under the global mean integrated squared error over a collection of spatially inhomogeneous function classes. Numerical implementation and simulation results are also discussed.

**1. Introduction.** Variance function estimation in heteroscedastic nonparametric regression is important in many contexts. In addition to being of interest in its own right, variance function estimates are needed, for example, to construct confidence intervals/bands for the mean function and to compute weighted least squares estimates of the mean function. Relative to mean function estimation, the literature on variance function estimation is sparse. Hall and Carroll (1989) considered kernel estimators of the variance function based on the squared residuals from a rate optimal estimator of the mean function. Müller and Stadtmüller (1987, 1993) considered difference based kernel estimators of the variance function. Ruppert et al. (1997) and Fan and Yao (1998) estimated the variance function by using local polynomial smoothing of the squared residuals from an "optimal" estimator of the mean function. More recently, Wang, Brown, Cai and Levine (2008)

Received May 2007; revised May 2007.
[1]Supported in part by NSF Grant DMS-06-04954.

derived the minimax rate of convergence for variance function estimation and constructed minimax rate optimal kernel estimators. Brown and Levine (2007) proposed a class of difference-based kernel estimators and established asymptotic normality.

So far the attention has been mainly focused on nonadaptive estimation of the variance function, that is, the smoothness of the variance function is assumed to be known and the estimators depend on the smoothness. In practice, however, the smoothness of the underlying functions is nearly always unknown. It is thus important to construct estimators that automatically adapt to the smoothness of the mean and variance functions. This is the goal of the present paper. Specifically, we propose a wavelet thresholding approach to adaptive variance function estimation in the heteroscedastic nonparametric regression model

$$(1) \qquad y_i = f(x_i) + V^{1/2}(x_i) z_i, \qquad i = 1, \ldots, n,$$

where $x_i = i/n$ and $z_i$ are independent and identically distributed with zero mean and unit variance. Here $n = 2^J$ for some positive integer $J$. The primary object of interest is the variance function $V(x)$. The estimation accuracy is measured both globally by the mean integrated squared error (MISE)

$$(2) \qquad R(\widehat{V}, V) = E\|\widehat{V} - V\|_2^2$$

and locally by the mean squared error (MSE) at a given point $x_* \in (0, 1)$

$$(3) \qquad R(\widehat{V}(x_*), V(x_*)) = E(\widehat{V}(x_*) - V(x_*))^2.$$

It is well known that when the mean function is sufficiently smooth, it has no first order effect on the quality of estimation for the variance function $V$; that is, one can estimate $V$ with the same asymptotic risk as if $f$ were known. See, for example, Ruppert et al. (1997) and Fan and Yao (1998). On the other hand, when $f$ is not smooth, the difficulty in estimating $V$ can be completely driven by the degree of smoothness of the mean $f$. How the smoothness of the unknown mean function influences the rate of convergence of the variance estimator can be characterized explicitly. Wang et al. (2008) showed that the minimax rate of convergence under both the pointwise MSE and global MISE is

$$(4) \qquad \max\{n^{-4\alpha}, n^{-2\beta/(2\beta+1)}\}$$

if $f$ has $\alpha$ derivatives and $V$ has $\beta$ derivatives.

The goal of the present paper is to estimate the variance function adaptively without assuming the degree of smoothness for either the mean function $f$ or variance function $V$. We introduce a wavelet thresholding procedure which applies wavelet thresholding to the squared first-order differences of the observations in (1). The procedure has two main steps. The first step is taking the square of the first-order differences of the observations $y_i$. This step turns the problem of variance function estimation under the model (1) into a more conventional regression problem of estimating the mean function. Another motivation for taking the differences is to eliminate the effect of the mean function $f$. The second step is to apply a wavelet thresholding procedure to the squared differences.

The procedure enjoys a high degree of adaptivity and spatial adaptivity in terms of the rate of convergence both for global and local estimation. More specifically, under the global risk measure (2), it adaptively achieves within a logarithmic factor of the minimax risk over a wide range of function classes which contain spatially inhomogeneous functions that may have, for example, jump discontinuities and high frequency oscillations. The estimator also optimally adapts to the local smoothness of the underlying function. As a special case, it is shown that the variance function estimator adaptively achieves the rate of convergence

$$(5) \qquad \max\left\{n^{-4\alpha}, \left(\frac{\log n}{n}\right)^{2\beta/(1+2\beta)}\right\}$$

under both the pointwise MSE and global MISE, if $f$ has $\alpha$ derivatives and $V$ has $\beta$ derivatives. Furthermore, it is shown that the extra logarithmic factor in the adaptive rate of convergence in (5) is necessary under the pointwise MSE and the estimator is thus optimally locally adaptive.

The wavelet estimator of the variance function is data-driven and easily implementable. We implement the procedure in S-Plus and R and carry out a simulation study to investigate the numerical performance of the estimator. Simulation results show that the MISE mostly depends on the structure of the underlying variance function and the effect of the mean function is not significant. In addition, we also compare the performance of the wavelet estimator with that of a kernel estimator whose bandwidth is chosen by cross validation. The numerical results show that the wavelet estimator uniformly outperforms the kernel estimator.

The paper is organized as follows. After Section 2.1 in which basic notation and definitions are summarized, the wavelet thresholding procedure is introduced in Sections 2.2 and 2.3. Sections 3 and 4 investigate the theoretical properties of the estimator. In particular, Section 4.1 derives a rate-sharp lower bound for the adaptive rate of convergence under the pointwise squared error loss. The lower and upper bounds together show that the estimator is optimally adaptive under the pointwise loss. Section 5 discusses implementation of the procedure and presents the numerical results. The proofs are given in Section 6.

**2. Wavelet procedure for variance function estimation.** In this section we introduce a wavelet thresholding procedure for estimating the variance function $V$ under the heteroscedastic regression model (1). We begin with notation and definitions of wavelets and a brief introduction to wavelet thresholding for estimating the mean function in the standard Gaussian regression setting and then give a detailed description of our wavelet procedure for variance function estimation.

2.1. *Wavelet thresholding for Gaussian regression.* We work with an orthonormal wavelet basis generated by dilation and translation of a compactly supported mother wavelet $\psi$ and a father wavelet $\phi$ with $\int \phi = 1$. A wavelet $\psi$ is called *r-regular* if $\psi$ has $r$ vanishing moments and $r$ continuous derivatives. A special family of compactly supported wavelets is the so-called Coiflets, constructed by Daubechies (1992), which can have arbitrary number of vanishing moments for both $\phi$ and $\psi$. Denote by $W(D)$ the collection of Coiflets $\{\phi, \psi\}$ of order $D$. So if $\{\phi, \psi\} \in W(D)$, then $\phi$ and $\psi$ are compactly supported and satisfy $\int x^i \phi(x)\,dx = 0$ for $i = 1, \ldots, D-1$; and $\int x^i \psi(x)\,dx = 0$ for $i = 0, \ldots, D-1$.

For simplicity in exposition, in the present paper we use periodized wavelet bases on $[0, 1]$. Let

$$\phi_{j,k}^p(x) = \sum_{l=-\infty}^{\infty} \phi_{j,k}(x - l),$$

$$\psi_{j,k}^p(x) = \sum_{l=-\infty}^{\infty} \psi_{j,k}(x - l) \qquad \text{for } t \in [0, 1],$$

where $\phi_{j,k}(x) = 2^{j/2}\phi(2^j x - k)$ and $\psi_{j,k}(x) = 2^{j/2}\psi(2^j x - k)$. The collection $\{\phi_{j_0,k}^p, k = 1, \ldots, 2^{j_0}; \psi_{j,k}^p, j \geq j_0 \geq 0, k = 1, \ldots, 2^j\}$ is then an orthonormal basis of $L^2[0,1]$, provided the primary resolution level $j_0$ is large enough to ensure that the support of the scaling functions and wavelets at level $j_0$ is not the whole of $[0, 1]$. The superscript "$p$" will be suppressed from the notation for convenience. An orthonormal wavelet basis has an associated orthogonal Discrete Wavelet Transform (DWT) which transforms sampled data into the wavelet coefficients. See Daubechies (1992) and Strang (1992) for further details about the wavelets and discrete wavelet transform. A square-integrable function $g$ on $[0, 1]$ can be expanded into a wavelet series:

(6) $$g(x) = \sum_{k=1}^{2^{j_0}} \xi_{j_0,k} \phi_{j_0,k}(x) + \sum_{j=j_0}^{\infty} \sum_{k=1}^{2^j} \theta_{j,k} \psi_{j,k}(x),$$

where $\xi_{j,k} = \langle g, \phi_{j,k} \rangle, \theta_{j,k} = \langle g, \psi_{j,k} \rangle$ are the wavelet coefficients of $g$.

Wavelet thresholding methods have been well developed for nonparametric function estimation, especially for estimating the mean function in the setting of homoscedastic Gaussian noise where one observes

$$(7) \qquad y_i = f\left(\frac{i}{n}\right) + \sigma z_i, \qquad z_i \stackrel{\text{i.i.d.}}{\sim} N(0,1), \qquad i = 1, \ldots, n.$$

One of the best known wavelet thresholding procedures is Donoho–Johnstone's VisuShrink [Donoho and Johnstone (1994) and Donoho et al. (1995)]. A typical wavelet thresholding procedure has three steps:

1. Transform the noisy data via the discrete wavelet transform.
2. Threshold the empirical wavelet coefficients by "killing" coefficients of small magnitude and keeping the large coefficients.
3. Estimate function $f$ at the sample points by inverse discrete wavelet transform of the denoised wavelet coefficients.

Many wavelet procedures are adaptive and easy to implement. We shall develop in Section 2.2 a wavelet procedure for variance function estimation where the noise is heteroscedastic and non-Gaussian.

2.2. *Adaptive wavelet procedure for estimating the variance function.* We now give a detailed description of our wavelet thresholding procedure for variance function estimation. The procedure begins with taking squared difference of the observations and then applies wavelet thresholding to obtain an estimator of the variance function.

Set $D_i = \frac{1}{\sqrt{2}}(y_{2i-1} - y_{2i})$ for $i = 1, 2, \ldots, n/2$. Then one can write

$$\begin{aligned}
D_i &= \frac{1}{\sqrt{2}}(f(x_{2i-1}) - f(x_{2i}) + V^{1/2}(x_{2i-1})z_{2i-1} - V^{1/2}(x_{2i})z_{2i}) \\
&= \frac{1}{\sqrt{2}}\delta_i + V_i^{1/2}\varepsilon_i,
\end{aligned} \qquad (8)$$

where $\delta_i = f(x_{2i-1}) - f(x_{2i})$, $V_i^{1/2} = \sqrt{\frac{1}{2}(V(x_{2i-1}) + V(x_{2i}))}$ and

$$(9) \qquad \varepsilon_i = (V(x_{2i-1}) + V(x_{2i}))^{-1/2}(V^{1/2}(x_{2i-1})z_{2i-1} - V^{1/2}(x_{2i})z_{2i})$$

has zero mean and unit variance. Then $D_i^2$ can be written as

$$(10) \qquad D_i^2 = V_i + \tfrac{1}{2}\delta_i^2 + \sqrt{2}V_i^{1/2}\delta_i\varepsilon_i + V_i(\varepsilon_i^2 - 1).$$

In the above expression $V_i$ is what we wish to estimate, $\frac{1}{2}\delta_i^2$ is a bias term caused by the mean function $f$, and $\sqrt{2}V_i^{1/2}\delta_i\varepsilon_i + V_i(\varepsilon_i^2 - 1)$ is viewed as the noise term. By taking squared differences, we have turned the problem of estimating the variance function into the problem of estimating the mean function similar to the conventional Gaussian regression model (7). The differences are of course that the noise is non-Gaussian and heteroscedastic and that there are additional unknown deterministic errors. In principle, virtually any good nonparametric regression procedure for estimating the mean function can then be applied. In this paper we shall use a wavelet estimator for its spatial adaptivity, asymptotic optimality, and computational efficiency.

We now construct a wavelet thresholding estimator $\widehat{V}$ based on the squared differences $D_i^2$. Although the procedure is more complicated, the basic idea is similar to the one behind the VisuShrink estimator for homoscedastic Gaussian regression described at the end of Section 2.1. We first apply the discrete wavelet transform to $\tilde{D} = \sqrt{2/n}(D_1^2, \ldots, D_{n/2}^2)'$. Let $d = W \cdot \tilde{D}$ be the empirical wavelet coefficients, where $W$ is the discrete wavelet transformation matrix. Then $d$ can be written as

$$(11) \quad d = (\tilde{d}_{j_0,1}, \ldots, \tilde{d}_{j_0,2^{j_0}}, d_{j_0,1}, \ldots, d_{j_0,2^{j_0}}, \ldots, d_{J-2,1}, \ldots, d_{J-2,2^{J-2}})',$$

where $\tilde{d}_{j_0,k}$ are the gross structure terms at the lowest resolution level, and $d_{j,k}$ ($j = j_0, \ldots, J-1, k = 1, \ldots, 2^j$) are empirical wavelet coefficients at level $j$ which represent fine structure at scale $2^j$. For convenience, we use $(j,k)$ to denote the number $2^j + k$. Then the empirical wavelet coefficients can be written as

$$d_{j,k} = \tau_{j,k} + z_{j,k}$$

where $z_{j,k}$ denotes the transformed noise part, that is,

$$z_{j,k} = \sqrt{\frac{2}{n}} \sum_i W_{(j,k),i}(\sqrt{2}V_i^{1/2}\delta_i\varepsilon_i + V_i(\varepsilon_i^2 - 1))$$

and

$$\tau_{j,k} = \theta_{j,k} + \sum_i W_{(j,k),i}\sqrt{\frac{1}{2n}}\delta_i^2 + \gamma_{j,k}.$$

Here $\theta_{j,k}$ is the true wavelet coefficients of $V(x)$, that is, $\theta_{j,k} = \langle V, \psi_{j,k} \rangle$, and $\gamma_{j,k}$ is the difference between $\theta_{j,k}$ and the discrete wavelet coefficient of $V_i$,

$$\gamma_{j,k} = \sum_i W_{(j,k),i}V_i - \theta_{j,k}.$$

We shall see that the approximation error $\gamma_{j,k}$ is negligible.

For the gross structure terms at the lowest resolution level, similarly, we can write

$$\tilde{d}_{j_0,k} = \tilde{\tau}_{j_0,k} + \tilde{z}_{j_0,k},$$

where

$$\widetilde{\tau}_{j_0,k} = \sum_i W_{(j_0,k),i}\sqrt{\frac{1}{2n}}\delta_i^2 + \xi_{j_0,k} + \widetilde{\gamma}_{j_0,k}, \tag{12}$$

$$\widetilde{z}_{j_0,k} = \sqrt{\frac{2}{n}}\sum_i W_{(j_0,k),i}(\sqrt{2}V_i^{1/2}\delta_i\varepsilon_i + V_i(\varepsilon_i^2-1)), \tag{13}$$

with $\xi_{j_0,k} = \langle V, \phi_{j_0,k}\rangle$ and $\widetilde{\gamma}_{j_0,k} = \sum_i W_{(j_0,k),i}V_i - \xi_{j_0,k}$.

Note that the squared differences $D_i^2$ are independent and the variance $\sigma_{j,k}^2$ of the empirical wavelet coefficients $d_{j,k}$ can be calculated as follows:

$$\sigma_{j,k}^2 \equiv \mathrm{Var}(d_{j,k}) = \frac{2}{n}\sum_i^{n/2} W_{(j,k),i}^2 \mathrm{Var}(D_i^2). \tag{14}$$

We shall use this formula to construct an estimator of $\sigma_{j,k}^2$ and then use it for choosing the threshold.

For any $y$ and $t \geq 0$, define the soft thresholding function $\eta_t(y) = \mathrm{sgn}(y)(|y|-t)_+$. Let $J_1$ be the largest integer satisfying $2^{J_1} \leq J^{-3}2^J$. Then the $\theta_{j,k}$ are estimated by

$$\hat{\theta}_{j,k} = \begin{cases} \eta_{\lambda_{j,k}}(d_{j,k}), & \text{if } j_0 \leq j \leq J_1, \\ 0, & \text{otherwise}, \end{cases} \tag{15}$$

where

$$\lambda_{j,k} = \hat{\sigma}_{j,k}\sqrt{2\log(n/2)} \tag{16}$$

is the threshold level, and $\hat{\sigma}_{j,k}$ is an estimate of the standard deviation $\sigma_{j,k}$. We shall discuss estimation of $\sigma_{j,k}$ in Section 2.3.

In wavelet regression for estimating the mean function, the coefficients of the father wavelets $\phi_{j_0,k}$ at the lowest resolution level are conventionally estimated by the corresponding empirical coefficients. Since there are only a small fixed number of coefficients, they would not affect the rate of convergence and the numerical results show that the wavelet estimators perform well in general. But in the setting of the present paper for estimating the variance function it turns out that using the empirical coefficients directly, although not affecting rate of convergence, often does not yield good numerical performance. We therefore estimate these coefficients by a shrinkage estimator given in Berger (1976). The estimator depends on the covariance matrix of the empirical coefficients $\widetilde{d}_{j_0} = (\widetilde{d}_{j_0,1},\ldots,\widetilde{d}_{j_0,2^{j_0}})'$ which is a function of the means and thus unknown. We shall use an estimated covariance matrix. Specifically, the estimator $\hat{\xi}_{j_0}$ of $\xi_{j_0}$ is given by

$$\hat{\xi}_{j_0} = \left(I - \frac{\min\{\widetilde{d}_{j_0}'\hat{\Sigma}^{-1}\widetilde{d}_{j_0}, 2^{j_0}-2\}\hat{\Sigma}^{-1}}{\widetilde{d}_{j_0}'\hat{\Sigma}^{-1}\hat{\Sigma}^{-1}\widetilde{d}_{j_0}}\right)\widetilde{d}_{j_0}, \tag{17}$$

where $\hat{\xi}_{j_0} = (\hat{\xi}_{j_0,1}, \hat{\xi}_{j_0,2}, \ldots, \hat{\xi}_{j_0,2^{j_0}})'$ is the estimator and $\hat{\Sigma}$ is the estimated covariance matrix of $\tilde{d}_{j_0}$. In our problem, we set

$$\hat{\Sigma} = \frac{2}{n} W_{j_0} \widehat{V}_D W'_{j_0}$$

where $W_{j_0}$ is the father wavelets part of the discrete wavelet transform matrix $W$. That is, $W_{j_0}$ is a $2^{j_0} \times \frac{n}{2}$ matrix and in our setting $W_{j_0}$ consists of the first $2^{j_0}$ rows of $W$. $\widehat{V}_D$ is a diagonal matrix given by

$$\widehat{V}_D = \text{Diag}\{\widehat{\text{Var}(D_1^2)}, \widehat{\text{Var}(D_2^2)}, \ldots, \widehat{\text{Var}(D_{n/2}^2)}\}$$

with $\widehat{\text{Var}(D_i^2)}$ given in equation (20) in Section 2.3.

With $\hat{\theta}_{j,k}$ given in (15) and $\hat{\xi}_{j_0,k}$ in (17), the estimator of the variance function $V$ is defined by

$$(18) \qquad \widehat{V}_e(x) = \sum_{k=1}^{2^{j_0}} \hat{\xi}_{j_0,k} \phi_{j_0,k}(x) + \sum_{j=j_0}^{J_1} \sum_{k=1}^{2^j} \hat{\theta}_{j,k} \psi_{j,k}(x).$$

So far we have only used half of the differences. Similarly we can apply the same procedure to the other half of differences, $D'_i = \frac{1}{\sqrt{2}}(y_{2i} - y_{2i+1})$, and obtain another wavelet estimator $\widehat{V}_o(x)$. The final estimator of the variance function $V$ is then given by

$$(19) \qquad \widehat{V}(x) = \tfrac{1}{2}(\widehat{V}_e(x) + \widehat{V}_o(x)).$$

The variance function at the sample points $\underline{V} = \{V(\frac{2i}{n}) : i = 1, \ldots, n/2\}$ can be estimated by the inverse transform of the denoised wavelet coefficients: $\underline{\widehat{V}} = W^{-1} \cdot (\frac{n}{2})^{1/2} \widehat{\Theta}$.

REMARK. In principle other thresholding techniques such as BlockJS [Cai (1999)], NeighBlock [Cai and Silverman (2001)] and EbayesThresh [Johnstone and Silverman (2005)] can also be adopted for estimating the wavelet coefficients here. However in the setting of the present paper the empirical coefficients are non-Gaussian, heteroscedastic and correlated and this makes the analysis of the properties of the resulting estimators very challenging.

2.3. *Estimate of $\sigma_{j,k}^2$*. Our wavelet estimator (19) of the variance function $V$ requires an estimate of variance $\sigma_{j,k}^2$ of the empirical wavelet coefficients. To make the wavelet estimator $\widehat{V}$ perform well asymptotically, we need a positively biased estimate of $\sigma_{j,k}^2$. That is, the estimate $\hat{\sigma}_{j,k}^2$ is greater than or equal to $\sigma_{j,k}^2$ with large probability. This can be seen in the proof of our theoretical results. On the other hand, the estimate $\hat{\sigma}_{j,k}^2$ should be of the same order as $\sigma_{j,k}^2$. Combining the two requirements together we need an estimate $\hat{\sigma}_{j,k}^2$ such that $P(\bigcap_{j,k} \sigma_{j,k}^2 \leq \hat{\sigma}_{j,k}^2 \leq C\sigma_{j,k}^2) \to 1$ for some constant $C > 1$. We shall construct such an estimate here.

It is clear from equation (14) that an estimate of the variance of $D_i^2$, $i = 1, 2, \ldots, n/2$, yields an estimate of $\sigma_{j,k}^2$. It follows from equation (10) that we need an estimate of $E(e_i^2)$, where $e_i = \sqrt{2}V_i^{1/2}\delta_i\varepsilon_i + V_i(\varepsilon_i^2 - 1)$. Note that $D_i^2 = V_i + \frac{1}{2}\delta_i^2 + e_i$ where $\delta_i$ is "small" and $V_i$ is "smooth." We shall estimate $\text{Var}(D_i^2) = E(e_i^2)$ by the average squared differences of $D_i^2$ over an subinterval. Specifically, we define the estimator of $\text{Var}(D_i^2)$ as follows.

Let $\Delta_i = D_{2i-1}^2 - D_{2i}^2$ for $i = 1, 2, \ldots, n/4$. Fix $0 < r < 1$ and divide the indices $1, 2, \ldots, n/2$ into nonoverlapping blocks of length $[(\frac{n}{2})^r]$. We estimate $\text{Var}(D_i^2) = E(e_i^2)$ in each block by the same value. Let $K$ be the total number of blocks and $B_k$ be the set of indices in the $k$th block. For $1 \leq k \leq K$, let

$$(20) \quad \widehat{\text{Var}(D_i^2)} \equiv \hat{\sigma}_k^2 = \frac{2}{(n/2)^r(2 - 1/\log n)} \sum_{2t \in B_k} \Delta_t^2 \qquad \text{for all } i \in B_k.$$

Lemma 6 in Section 6 shows that this estimate has the desired property for any $0 < r < 1$. With this estimator of $\text{Var}(D_i^2)$, we estimate $\sigma_{j,k}^2$ by

$$(21) \qquad \hat{\sigma}_{j,k}^2 = \frac{2}{n} \sum_{i}^{n/2} W_{(j,k),i}^2 \widehat{\text{Var}(D_i^2)}.$$

We shall use $\hat{\sigma}_{j,k}$ in the threshold $\lambda_{j,k}$ in (16) and construct the wavelet estimator of $V$ as in (18) and (19).

**3. Global adaptivity of the wavelet procedure.** We consider in this section the theoretical properties of the variance function estimator $\widehat{V}$ given in (19) under the global MISE (2). The local adaptivity of the estimator under pointwise MSE (3) is treated in Section 4. These results show that the variance function estimator (19) is nearly optimally adaptive and spatially adaptive over a wide range of function spaces for both the mean and variance functions.

3.1. *Inhomogeneous function class $\mathcal{H}$.* To demonstrate the global adaptivity of the variance function estimator $\widehat{V}$, we consider a family of large function classes which contain spatially inhomogeneous functions that may have, for example, jump discontinuities and high frequency oscillations. These function classes are different from the more traditional smoothness classes. Functions in these classes can be viewed as the superposition of smooth functions with irregular perturbations. These and other similar function classes have been used in Hall, Kerkyacharian and Picard (1998, 1999) and Cai (2002) in the study of wavelet block thresholding estimators.

DEFINITION 1. Let $\mathcal{H} = \mathcal{H}(\alpha_1, \alpha, \gamma, M_1, M_2, M_3, D, v)$, where $0 \leq \alpha_1 < \alpha \leq D$, $\gamma > 0$, and $M_1, M_2, M_3, v \geq 0$, denote the class of functions $f$ such that for any $j \geq j_0 > 0$ there exists a set of integers $A_j$ with $\text{card}(A_j) \leq M_3 2^{j\gamma}$ for which the following are true:

- For each $k \in A_j$, there exist constants $a_0 = f(2^{-j}k), a_1, \ldots, a_{D-1}$ such that for all $x \in [2^{-j}k, 2^{-j}(k+v)]$, $|f(x) - \sum_{m=0}^{D-1} a_m(x - 2^{-j}k)^m| \leq M_1 2^{-j\alpha_1}$.
- For each $k \notin A_j$, there exist constants $a_0 = f(2^{-j}k), a_1, \ldots, a_{D-1}$ such that for all $x \in [2^{-j}k, 2^{-j}(k+v)]$, $|f(x) - \sum_{m=0}^{D-1} a_m(x - 2^{-j}k)^m| \leq M_2 2^{-j\alpha}$.

A function $f \in \mathcal{H}(\alpha_1, \alpha, \gamma, M_1, M_2, M_3, D, v)$ can be regarded as the superposition of a regular function $f_s$ and an irregular perturbation $\tau$: $f = f_s + \tau$. The perturbation $\tau$ can be, for example, jump discontinuities or high frequency oscillations such as chirp and Doppler of the form: $\tau(x) = \sum_{k=1}^{K} a_k(x - x_k)^{\beta_k} \cos(x - x_k)^{-\gamma_k}$. The smooth function $f_s$ belongs to the conventional Besov class $B_{\infty\infty}^{\alpha}(M_2)$. Roughly speaking, a Besov space $B_{p,q}^{\alpha}$ contains functions having $\alpha$ bounded derivatives in $L^p$ space, the parameter $q$ gives a finer gradation of smoothness. See Triebel (1983) and Meyer (1992) for more details on Besov spaces.

Intuitively, the intervals with indices in $A_j$ are "bad" intervals which contain less smooth parts of the function. The number of the "bad" intervals is controlled by $M_3$ and $\gamma$ so that the irregular parts do not overwhelm the fundamental structure of the function. It is easy to see that $\mathcal{H}(\alpha_1, \alpha, \gamma, M_1, M_2, M_3, D, v)$ contains the Besov class $B_{\infty\infty}^{\alpha}(M_2)$ as a subset for any given $\alpha_1$, $\gamma$, $M_1$, $M_3$, $D$ and $v$. See Hall, Kerkyacharian and Picard (1999) for further discussions on the function classes $\mathcal{H}$.

3.2. *Global adaptivity.* The minimax rate of convergence for estimating the variance function $V$ over the traditional Lipschitz balls was derived in Wang et al. (2008). Define the Lipschitz ball $\Lambda^{\alpha}(M)$ in the usual way as

$$\Lambda^{\alpha}(M) = \{g : |g^{(k)}(x)| \leq M, \text{ and } |g^{(\lfloor \alpha \rfloor)}(x) - g^{(\lfloor \alpha \rfloor)}(y)| \leq M|x-y|^{\alpha'}$$
$$\text{for all } 0 \leq x, y \leq 1, k = 0, \ldots, \lfloor \alpha \rfloor - 1\}$$

where $\lfloor \alpha \rfloor$ is the largest integer less than $\alpha$ and $\alpha' = \alpha - \lfloor \alpha \rfloor$. Wang et al. (2008) showed that the minimax risks for estimating $V$ over $f \in \Lambda^{\alpha}(M_f)$ and $V \in \Lambda^{\beta}(M_V)$ under both the global MISE and the MSE at a fixed point $x_* \in (0, 1)$ satisfy

$$\inf_{\widehat{V}} \sup_{f \in \Lambda^{\alpha}(M_f), V \in \Lambda^{\beta}(M_V)} E\|\widehat{V} - V\|_2^2$$
(22)
$$\asymp \inf_{\widehat{V}} \sup_{f \in \Lambda^{\alpha}(M_f), V \in \Lambda^{\beta}(M_V)} E(\widehat{V}(x_*) - V(x_*))^2$$
$$\asymp \max\{n^{-4\alpha}, n^{-2\beta/(\beta+1)}\}.$$

We now consider the function class $\mathcal{H}(\alpha_1, \alpha, \gamma, M_1, M_2, M_3, D, v)$ defined in Section 3.1. Let $\mathcal{H}_f(\alpha) = \mathcal{H}(\alpha_1, \alpha, \gamma_f, M_{f1}, M_{f2}, M_{f3}, D_f, v_f)$ and $\mathcal{H}_V(\beta) = \mathcal{H}(\beta_1, \beta, \gamma_V, M_{V1}, M_{V2}, M_{V3}, D_V, v_V)$. Since $\mathcal{H}(\alpha_1, \alpha, \gamma, M_1, M_2, M_3, D, v)$ contains the Lipschitz ball $\Lambda^\alpha(M_2)$ as a subset for any given $\alpha_1$, $\gamma$, $M_1$, $M_3$, $D$, and $v$, a minimax lower bound for estimating $V$ over $f \in \mathcal{H}_f(\alpha)$ and $V \in \mathcal{H}_V(\beta)$ follows directly from (22):

$$(23) \quad \inf_{\widehat{V}} \sup_{f \in \mathcal{H}_f(\alpha), V \in \mathcal{H}_V(\beta)} E\|\widehat{V} - V\|_2^2 \geq C \cdot \max\{n^{-4\alpha}, n^{-2\beta/(1+2\beta)}\}.$$

The following theorem shows that the variance function estimator $\widehat{V}$ is adaptive over a range of the function classes $\mathcal{H}$. We shall assume that the error $z_i$ in the regression model (1) satisfies the property that the moment generating function of $z_i^2$, $G(x) = E(e^{xz_i^2})$, exists when $|x| < \rho$ for some constant $\rho > 0$. This condition implies that the moment generating function of $\varepsilon_i$ in (9), $G_\varepsilon(x) = E(e^{x\varepsilon_i})$, exists in a neighborhood of 0.

THEOREM 1. *Let $\{y_1, \ldots, y_n\}$ be given as in (1). Suppose the wavelets $\{\phi, \psi\} \in W(D)$ and the moment generating function of $z_i^2$ exists in a neighborhood of the origin. Suppose also $\gamma_f \leq 1 + 4\alpha_1 - 4\alpha$, and $\gamma_V \leq \frac{1+2\beta_1}{1+2\beta}$. Then the variance function estimator $\widehat{V}$ given in (19) satisfies that for some constant $C_0 > 0$ and all $0 < \beta \leq D$*

$$(24) \quad \sup_{f \in \mathcal{H}_f(\alpha), V \in \mathcal{H}_V(\beta)} E\|\widehat{V} - V\|_2^2 \leq C_0 \cdot \max\left\{n^{-4\alpha}, \left(\frac{\log n}{n}\right)^{2\beta/(1+2\beta)}\right\}.$$

REMARK. The use of Coiflets in Theorem 1 is purely for technical reasons. If the following mild local Lipschitz condition is imposed on functions in $\mathcal{H}$ in regions where the functions are relatively smooth, then the Coiflets are not needed. Local adaptivity result given in the next section does not require the use of Coiflets and our simulation shows no particular advantages of using Coiflets in the finite sample case.

(i) If $\alpha > 1 \geq \alpha_1$, then for $k \notin A_j, |f(x) - f(2^{-j}k)| \leq M_4 2^{-j}$, for $x \in [2^{-j}k, 2^{-j}(k+v)]$.

(ii) If $\alpha > \alpha_1 > 1$, then $|f(x) - f(2^{-j}k)| \leq M_4 2^{-j}$, for $x \in [2^{-j}k, 2^{-j}(k+v)]$.

Comparing (24) with the minimax lower bound given in (23), the estimator $\widehat{V}$ is adaptively within a logarithmic factor of the minimax risk under global MISE. Thus, the variance function estimator $\widehat{V}$, without knowing the a priori degree or amount of smoothness of the underlying mean and variance functions, achieves within a logarithmic factor of the true optimal convergence rate that one could achieve by knowing the regularity.

For adaptive estimation of $V$ over the traditional Lipschitz balls, the following is a direct consequence of Theorem 1.

COROLLARY 1. *Let $\{y_1, \ldots, y_n\}$ be given as in (1). Suppose the wavelet $\psi$ is $r$-regular and the moment generating function of $z_i^2$ exists in a neighborhood of the origin. Then the variance function estimator $\widehat{V}$ given in (19) satisfies that for some constant $C_0 > 0$ and all $0 < \beta \leq r$*

$$(25) \quad \sup_{f \in \Lambda^\alpha(M_f), V \in \Lambda^\beta(M_V)} E\|\widehat{V} - V\|_2^2 \leq C_0 \cdot \max\left\{n^{-4\alpha}, \left(\frac{\log n}{n}\right)^{2\beta/(1+2\beta)}\right\}.$$

**4. Local adaptivity.** For functions of spatial inhomogeneity, the local smoothness of the functions varies significantly from point to point and global risk measures such as (2) cannot wholly reflect the performance of an estimator locally. The local risk measure (3) is more appropriate for measuring the spatial adaptivity, where $x_* \in (0,1)$ is any fixed point of interest.

Define the local Lipschitz class $\Lambda^\alpha(M, x_*, \delta)$ by

$$\Lambda^\alpha(M, x_*, \delta) = \{g : |g^{(\lfloor\alpha\rfloor)}(x) - g^{(\lfloor\alpha\rfloor)}(x_*)| \leq M|x - x_*|^{\alpha'},$$
$$x \in (x_* - \delta, x_* + \delta)\}$$

where $\lfloor\alpha\rfloor$ is the largest integer less than $\alpha$ and $\alpha' = \alpha - \lfloor\alpha\rfloor$.

THEOREM 2. *Let $\{y_1, \ldots, y_n\}$ be given as in (1). Suppose the wavelet $\psi$ is $r$-regular and the moment generating function of $z_i^2$ exists in a neighborhood of the origin. Then the variance function estimator $\widehat{V}$ given in (19) satisfies that for any fixed $x_* \in (0,1)$ there exists some constant $C_0 > 0$ such that for all $\alpha > 0$ and all $0 < \beta \leq r$*

$$(26) \quad \sup_{f \in \Lambda^\alpha(M_f, x_*, \delta_f), V \in \Lambda^\beta(M_V, x_*, \delta_V)} E(\widehat{V}(x_*) - V(x_*))^2$$
$$\leq C_0 \cdot \max\left\{n^{-4\alpha}, \left(\frac{\log n}{n}\right)^{2\beta/(1+2\beta)}\right\}.$$

Comparing (26) with the minimax rate of convergence given in (22), the estimator $\widehat{V}$ is simultaneously within a logarithmic factor of the minimax risk under the pointwise risk. We shall show in Section 4.1 that, under the pointwise risk, this logarithmic factor is unavoidable for adaptive estimation. It is the minimum penalty for not knowing the smoothness of the variance function $V$. Therefore the estimator $\widehat{V}$ is optimally adaptive under the pointwise loss.

4.1. *Lower bound for adaptive pointwise estimation.* We now turn to the lower bound for adaptive estimation of the variance function $V$ under the pointwise MSE. The sharp lower bound we derive below demonstrates that the cost of adaptation for variance function estimation behaves in a more complicated way than that for mean function estimation.

It is well known in estimating the mean function $f$ that it is possible to achieve complete adaptation for free under the global MISE in terms of the rate of convergence over a collection of function classes. That is, one can do as well when the degree of smoothness is unknown as one could do if the degree of smoothness is known. But for estimation at a point, one must pay a price for adaptation. The optimal rate of convergence for estimating the mean function $f$ at point over $\Lambda^\alpha(M_f)$ with $\alpha$ completely known is $n^{-2\alpha/(1+2\alpha)}$. In the setting of adaptive estimation of the mean function, Lepski (1990) and Brown and Low (1996) showed that one has to pay a price for adaptation of at least a logarithmic factor even when $\alpha$ is known to be one of two values. It is shown that the best achievable rate is $(\frac{\log n}{n})^{2\alpha/(1+2\alpha)}$, when the smoothness parameter $\alpha$ is unknown.

Here we consider adaptive estimation of the variance function $V$ at a point. The following lower bound characterizes the cost of adaptation for such a problem.

THEOREM 3. *Let $\alpha_0, \alpha_1 > 0$, $\beta_0 > \beta_1 > 0$ and $4\alpha_0 > \frac{2\beta_1}{1+2\beta_1}$. Under the regression model (1) with $z_i \stackrel{\text{i.i.d.}}{\sim} N(0,1)$, for any estimator $\widehat{V}$ and any fixed $x_* \in (0,1)$, if*

$$\overline{\lim_{n \to \infty}} \min\{n^{4\alpha_0}, n^{2\beta_0/(1+2\beta_0)}\}$$
(27)
$$\times \sup_{f \in \Lambda^{\alpha_0}(M_f), V \in \Lambda^{\beta_0}(M_V)} E(\widehat{V}(x_*) - V(x_*))^2 < \infty,$$

*then*

$$\underline{\lim_{n \to \infty}} \min\left\{n^{4\alpha_1}, \left(\frac{n}{\log n}\right)^{2\beta_1/(1+2\beta_1)}\right\}$$
(28)
$$\times \sup_{f \in \Lambda^{\alpha_1}(M_f), V \in \Lambda^{\beta_1}(M_V)} E(\widehat{V}(x_*) - V(x_*))^2 > 0.$$

The lower bound for adaptive estimation given in Theorem 3 is more complicated than the corresponding lower bound for estimating the mean function given in Lepski (1990) and Brown and Low (1996). Theorem 3 shows that, if an estimator is rate optimal for $f \in \Lambda^{\alpha_0}(M_f)$ and $V \in \Lambda^{\beta_0}(M_V)$, then one must pay a price of at least a logarithmic factor for $f \in \Lambda^{\alpha_1}(M_f)$ and $V \in \Lambda^{\beta_1}(M_V)$ if the mean function is smooth, that is, $4\alpha_1 \geq \frac{2\beta_1}{1+2\beta_1}$. On the other hand, if $4\alpha_1 < \frac{2\beta_1}{1+2\beta_1}$, then it is possible to achieve the exact minimax rate simultaneously both over $f \in \Lambda^{\alpha_0}(M_f)$, $V \in \Lambda^{\beta_0}(M_V)$ and $f \in \Lambda^{\alpha_1}(M_f)$, $V \in \Lambda^{\beta_1}(M_V)$. In contrast, for estimating the mean function at a point one must always pay a price of at least a logarithmic factor for not knowing the exact smoothness of the function.

Comparing the lower bound (28) with the upper bound (26) given in Theorem 2, it is clear that our wavelet estimator (19) is optimally adaptive under the pointwise risk. The lower and upper bounds together show the following. When the mean function is not smooth, that is, $4\alpha < \frac{2\beta}{1+2\beta}$, the minimax rate of convergence can be achieved adaptively. On the other hand, when the effect of the mean function is negligible, that is, $4\alpha \geq \frac{2\beta}{1+2\beta}$, the minimax rate of convergence cannot be achieved adaptively and one has to pay a minimum of a logarithmic factor as in the case of mean function estimation.

The proof of this theorem can be naturally divided into two parts. The first part

$$(29) \quad \lim_{n\to\infty} n^{4\alpha_1} \cdot \sup_{f\in\Lambda^{\alpha_1}(M_f), V\in\Lambda^{\beta_1}(M_V)} E(\widehat{V}(x_*) - V(x_*))^2 > 0$$

follows directly from the minimax lower bound given in Wang et al. (2008). We shall use a two-point constrained risk inequality to prove the second part,

$$(30) \quad \lim_{n\to\infty} \left(\frac{n}{\log n}\right)^{2\beta_1/(1+2\beta_1)} \cdot \sup_{f\in\Lambda^{\alpha_1}(M_f), V\in\Lambda^{\beta_1}(M_V)} E(\widehat{V}(x_*) - V(x_*))^2 > 0.$$

A detailed proof is given in Section 6.4.

**5. Numerical results.** The adaptive procedure for estimating the variance function introduced in Section 2.2 is easily implementable. We implement the procedure in S-Plus and R. In this section we will investigate the numerical performance of the estimator. The numerical study has three goals. The first is to investigate the effect of the mean function on the estimation of the variance function. Several different combinations of the mean and variance functions are used and the MSE of each case is given. The second goal is to study the effect of different choices of $r$ in (20) on the performance of the estimator. The simulation results indicate that the MISE of the estimator is not sensitive to the choice of $r$. Finally, we will make a comparison between the wavelet estimator and a kernel estimator with the bandwidth chosen by cross validation. For reasons of space, we only report here a summary of the numerical results. See Cai and Wang (2007) for more detailed and additional simulation results.

Four different variance functions were considered in the simulation study. They are Bumps and Doppler functions from Donoho and Johnstone (1994) and also the following two functions:

$$V_1(x) = \begin{cases} 3 - 30x, & \text{for } 0 \leq x \leq 0.1, \\ 20x - 1, & \text{for } 0.1 \leq x \leq 0.25, \\ 4 + (1 - 4x)18/19, & \text{for } 0.25 < x \leq 0.725, \\ 2.2 + 10(x - 0.725), & \text{for } 0.725 < x \leq 0.89, \\ 3.85 - 85(x - 0.89)/11, & \text{for } 0.89 < x \leq 1, \end{cases}$$

$$V_2(x) = 1 + 4(e^{-550(x-0.2)^2} + e^{-200(x-0.8)^2} + e^{-950(x-0.8)^2}).$$

These test functions are rescaled in the simulations to have the same $L_2$ norm.

We begin by considering the effect of the mean function on the estimation of the variance function. For each variance function $V(x)$, we use five different mean functions, the constant function $f(x) = 0$, the trigonometric function $f = \sin(20x)$, and Bumps, Blocks and Doppler functions from Donoho and Johnstone (1994). Different combinations of wavelets and sample size $n$ yield basically the same qualitative results. As an illustration, Table 1 reports the average squared errors over 500 replications with sample size $n = 4096$ using Daubechies compactly supported wavelet *Symmlet* 8. In this part, we use $r = 0.5$ in (21). Figure 1 provides a graphical comparison of the variance function estimators and the true functions in the case the mean function $f \equiv 0$.

It can be seen from Table 1 that in all these examples the MISEs mostly depend on the structure of the variance function. The effect of the mean function $f$ is not significant. For Bumps and Blocks, the spatial structure of the mean $f$ only affect a small number of wavelet coefficients, and the variance function estimator still performs well. But still, when $f$ is smooth, the estimator of the variance function $V$ is slightly more accurate. We can also see that the results here are not as good as the estimation of mean function under the standard homoscedastic Gaussian regression model. This is primarily due to the difficulty of the variance function estimation problem itself.

TABLE 1
*The average squared error over* 500 *replications with sample size* $n = 4096$

|          | $f(x) \equiv 0$ | $f(x) = \sin(20x)$ | Bumps  | Blocks | Doppler |
|----------|-----------------|---------------------|--------|--------|---------|
| $V_1(x)$ | 0.0817          | 0.0842              | 0.0825 | 0.0860 | 0.0837  |
| $V_2(x)$ | 0.0523          | 0.0553              | 0.0557 | 0.0563 | 0.0567  |
| Bumps    | 0.1949          | 0.2062              | 0.2146 | 0.2133 | 0.2060  |
| Doppler  | 0.4162          | 0.5037              | 0.4817 | 0.4888 | 0.4902  |

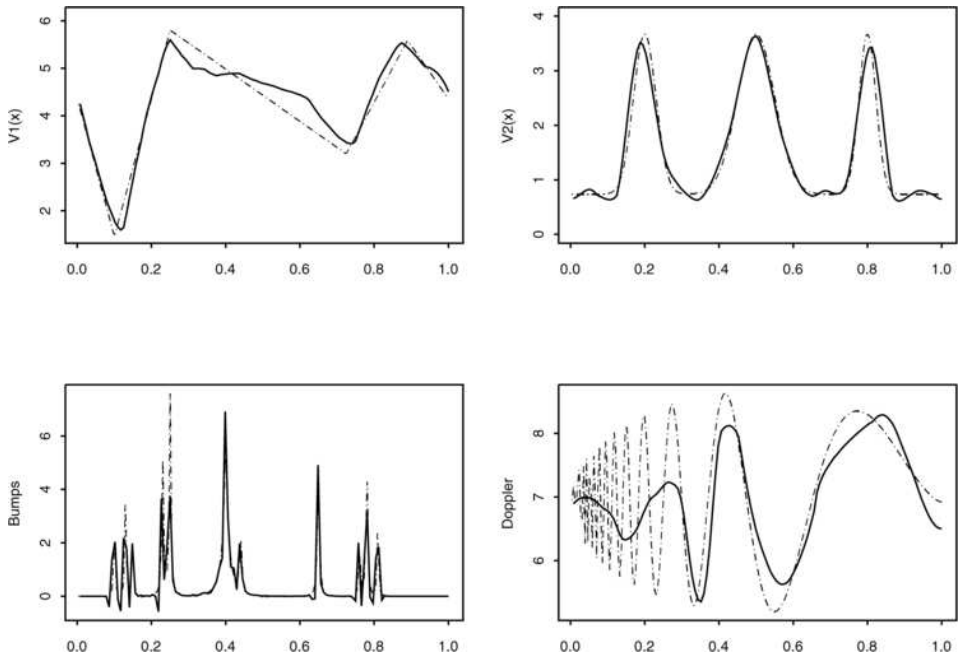

Fig. 1.  *Wavelet estimates (solid) and true variance functions (dotted).*

We now turn to the choice of $r$ in (21). Using the same setting as in the previous example, we apply our procedure for the four test functions with three different choices of $r$ in (21), $r = 0.2, 0.5$ and $0.8$, respectively. The mean function is chosen to be $f \equiv 0$. The average squared error over 500 replications are given in Table 2.

For each test function the MISEs are nearly identical for different choices of $r$. It is thus clear from Table 2 that the performance of the estimator is not sensitive to the choice of $r$. We suggest use $r = 0.5$ in practice.

After taking squared differences, the problem of estimating the variance function becomes the problem of estimating the mean function and virtually any good procedure for estimating the mean function can then be applied. We now compare the performance of our wavelet estimator with a kernel

TABLE 2
*The MISEs for different choices of r*

|  | $V_1(x)$ | $V_2(x)$ | Bumps | Doppler |
| --- | --- | --- | --- | --- |
| $r = 0.2$ | 0.0838 | 0.0581 | 0.1981 | 0.4852 |
| $r = 0.5$ | 0.0817 | 0.0523 | 0.1949 | 0.4162 |
| $r = 0.8$ | 0.0859 | 0.0532 | 0.2065 | 0.4335 |

ADAPTIVE VARIANCE FUNCTION ESTIMATION 17

Table 3
*Comparison of the MISEs for the wavelet and kernel estimators*

|         | $V_1(x)$ | $V_2(x)$ | Bumps  | Doppler |
|---------|----------|----------|--------|---------|
| Wavelet | 0.0817   | 0.0523   | 0.1949 | 0.4762  |
| Kernel  | 0.1208   | 0.0631   | 0.2296 | 0.5463  |

estimator whose bandwidth is chosen via cross-validation. Table 3 displays the average squared errors over 500 replications of the two estimators for the four variance functions with the mean function $f \equiv 0$.

The wavelet estimator outperforms the kernel estimator for all the variance functions. The MISEs of the kernel estimator are 14% to 47% higher than the corresponding wavelet estimator. Although the bandwidth of the kernel estimator is chosen adaptive via cross-validation, the spatial inhomogeneity of the variance functions limits the performance of any kernel method with a single bandwidth.

In summary, the simulation study shows that the effect of the mean function on the performance of the wavelet estimator is not significant. In this sense our wavelet procedure is robust against the mean function interference. The procedure is also not sensitive to the choice of $r$. In addition, the wavelet estimator uniformly outperforms the kernel estimator whose bandwidth is chosen by cross-validation.

**6. Proofs.** We begin by introducing and proving several technical lemmas in Section 6.1 that will be used in the proof of the main results. Throughout this section, we use $C$ (as well as $C_0$, $C_1$, etc.) to denote constants that may vary from place to place.

6.1. *Preparatory results.* Oracle inequality for the soft thresholding estimator was given in Donoho and Johnstone (1994) in the case when the noise is i.i.d. normal. In the present paper we need the following risk bound for the soft thresholding estimator without the normality assumption. This risk bound is useful in turning the analysis of the variance function estimator into the bias-variance trade-off calculation which is often used in more standard Gaussian nonparametric regression.

LEMMA 1. *Let $y = \theta + Z$, where $\theta$ is an unknown parameter and $Z$ is a random variable with $EZ = 0$. Then*
$$E(\eta(y,\lambda) - \theta)^2 \leq \theta^2 \wedge (4\lambda^2) + 2E(Z^2 I(|Z| > \lambda)).$$

PROOF. Note that
$$E(\eta(y,\lambda) - \theta)^2 \leq 2E(\eta(y,\lambda) - y)^2 + 2E(y - \theta)^2 \leq 2\lambda^2 + 2EZ^2$$
$$\leq 4\lambda^2 + 2E(Z^2 I(|Z| > \lambda)).$$

On the other hand,

$$\begin{aligned}E(\eta(y,\lambda)-\theta)^2 &= \theta^2 P(-\lambda-\theta \le Z \le \lambda-\theta) + E((Z-\lambda)^2 I(Z > \lambda-\theta)) \\ &\quad + E((Z+\lambda)^2 I(Z < -\lambda-\theta)) \\ &\le \theta^2 + E((Z-\lambda)^2 I(Z > \lambda)) + E((Z+\lambda)^2 I(Z < -\lambda)) \\ &\le \theta^2 + E(Z^2 I(|Z| > \lambda)). \quad \square\end{aligned}$$

The following lemma bounds the wavelet coefficients of the functions in $\mathcal{H}$.

LEMMA 2.  (i) Let $g \in \mathcal{H}(\alpha_1, \alpha, \gamma, M_1, M_2, M_3, D, v)$. Assume the wavelets $\{\phi, \varphi\} \in W(D)$ with $\operatorname{supp}(\phi) = \operatorname{supp}(\psi) \subset [0, v]$. Let $n = 2^J$, $\xi_{J,k} = \int g\phi_{J,k}$ and $\theta_{j,k} = \int g\psi_{j,k}$. Then

$$\begin{aligned}|\xi_{J,k} - n^{-1/2} g(k/n)| &\le M_1 \|\phi\|_1 n^{-(1/2+\alpha_1)} && \text{for all } k \in A_J; \\ |\xi_{J,k} - n^{-1/2} g(k/n)| &\le M_2 \|\phi\|_1 n^{-(1/2+\alpha)} && \text{for all } k \notin A_J; \\ |\theta_{j,k}| &\le M_1 \|\psi\|_1 2^{-j(1/2+\alpha_1)} && \text{for all } k \in A_j; \\ |\theta_{j,k}| &\le M_1 \|\psi\|_1 2^{-j(1/2+\alpha)} && \text{for all } k \notin A_j.\end{aligned}$$

(ii) *For all functions $g \in \Lambda^\alpha(M)$, the wavelet coefficients of $g$ satisfy* $|\theta_{j,k}| \le C 2^{-j(1/2+\alpha)}$ *where constant $C$ depends only on the wavelets, $\alpha$ and $M$ only.*

Lemma 2(ii) is a standard result; see for example, Daubechies (1992). For a proof of Lemma 2(i), see Hall, Kerkyacharian and Picard (1999) and Cai (2002). It follows from this lemma that

$$\sup_{g \in \Lambda^\beta(M)} \sum_{k=1}^n \left(\xi_{J,k} - n^{-1/2} g\left(\frac{k}{n}\right)\right)^2 \le C n^{-(2\beta \wedge 1)}. \tag{31}$$

The next lemma gives a large deviation result, which will be used to control the tail probability of the empirical wavelet coefficients.

LEMMA 3.  *Suppose $\varepsilon_i$, $i = 1, 2, \ldots$, are independent random variables with $E\varepsilon_i = 0$, $\operatorname{Var}(\varepsilon_i) = v_i \le v_0$ for all $i$. Moreover, suppose the moment generating function $M^i(x) \triangleq E(\exp(x\varepsilon_i))$ exists when $|x| < \rho$ for some $\rho > 0$ and all $i$. Let*

$$Z_m = \frac{1}{\sqrt{v_0}} \sum_{i=1}^m a_{mi} \varepsilon_i$$

with $\sum_{i=1}^{m} a_{mi}^2 = 1$ and $|a_{mi}| \leq c_0/\sqrt{m}$ for some constant $c_0$, then for $\lambda = o(m^{-1/4})$ and sufficiently large $m$

$$\frac{P(|Z_m| > \sigma_m \lambda)}{2(1 - \Phi(\lambda))} \leq \exp\left(C\frac{\lambda^3}{m^{1/4}}\right)(1 + O(m^{-1/4}))$$

where $\sigma_m^2 = \sum a_{mi}^2 v_i / v_0$ and $C > 0$ is a constant.

PROOF. Let $S_m = m^{1/4} Z_m$ and $M_m(x) = E(\exp(xS_m))$. Suppose $\mu_{ik}$ denote the $k$th moment of $\varepsilon_i$. Note that $|a_{mi}| \leq c_0 m^{-1/2}$ for any $m$ and $1 \leq i \leq m$, we have

$$\frac{\log(M_m(x))}{m^{1/2}} = \frac{1}{m^{1/2}} \sum_{i=1}^{m} \log M^i\left(x \frac{a_{mi} m^{1/4}}{\sqrt{v_0}}\right)$$

$$= \frac{1}{m^{1/2}} \sum_{i=1}^{m} \log\left(1 + \frac{m^{1/2} a_{mi}^2}{2v_0} x^2 + \frac{m^{3/4} a_{mi}^3}{6 v_0^{3/2}} \mu_{i3} x^3 + \cdots\right)$$

$$= \frac{1}{m^{1/2}} \sum_{i=1}^{m} \left(\frac{m^{1/2} a_{mi}^2}{2v_0} x^2 + m^{-3/4} x^3 \cdot \Delta_m(x)\right)$$

$$= \frac{x^2}{2v_0} + m^{-1/4} x^3 \cdot \Delta_m(x)$$

where $\Delta_m(x)$ is uniformly bounded for all $m$ when $x < \rho$. This means that $M_m(x)$ can be written in the form $M_m(x) = e^{m^{1/2}(x^2/2v_0)}(1 + O(m^{-1/4}))$ for $|x| < \rho$. It then follows from Theorem 1 of Hwang (1996) that for $\lambda = o(m^{-1/4})$ and sufficiently large $m$

$$\frac{P(|Z_m| > \sigma_m \lambda)}{2(1 - \Phi(\lambda))} \leq \exp\left(C\frac{\lambda^3}{m^{1/4}}\right)(1 + O(m^{-1/4})). \qquad \square$$

A special case of Lemma 3 is when $m \geq (\log n)^k$ for a positive integer $n$ and some $k > 2$ and $\lambda = \sqrt{2 \log n}$. In this case, we have

$$\frac{P(|Z_m| > \sigma_m \sqrt{2 \log n})}{2(1 - \Phi(\sqrt{2 \log n}))} \leq \exp\left(C\frac{(2 \log n)^{3/2}}{(\log n)^{k/4}}\right)\left(1 + O\left(\frac{\sqrt{2 \log n}}{(\log n)^{k/4}}\right)\right).$$

Since $k > 2$, $\exp\{(\log n)^{3/2 - k/4}\} = o(n^a)$ for any $a > 0$ as $n \to \infty$. Therefore

(32) $$P(|Z_m| > \sigma_m \sqrt{2 \log n}) \leq O\left(\frac{1}{n^{1-a}}\right)$$

for any $a > 0$.

The following two lemmas bounds the difference between the mean $\tau_{j,k}$ of the empirical wavelet coefficient $d_{j,k}$ and the true wavelet coefficient $\theta_{j,k}$, globally and individually.

LEMMA 4. *Using the notation in Section 2.2, and under the conditions of Theorem 1, we have*

$$\sup_{f\in\mathcal{H}_f(\alpha), V\in\mathcal{H}_V(\beta)} \left\{ \sum_k (\tilde{\tau}_{j_0,k} - \xi_{j_0,k})^2 + \sum_{j,k} (\tau_{j,k} - \theta_{j,k})^2 \right\}$$
$$= O(n^{-(1\wedge 4\alpha \wedge 2\beta \wedge (1+2\beta_1 - \gamma_V))}).$$

PROOF. Note that

$$\sum_k (\tilde{\tau}_{j_0,k} - \xi_{j_0,k})^2 + \sum_{j,k}(\tau_{j,k} - \theta_{j,k})^2 = \sum_{j,k}\left(\sum_i W_{(j,k),i}\sqrt{\frac{1}{2n}}\delta_i^2 + \gamma_{j,k}\right)^2$$

$$\leq 2\sum_{j,k}\left(\sum_i W_{(j,k),i}\sqrt{\frac{1}{2n}}\delta_i^2\right)^2 + 2\sum_{j,k}\gamma_{j,k}^2.$$

It follows from the isometry property of the orthogonal wavelet transform that

$$\sum_{j,k}\left(\sum_i W_{(j,k),i}\sqrt{\frac{1}{2n}}\delta_i^2\right)^2 = \frac{1}{2n}\sum_i \delta_i^4.$$

From the definition of the function class $\mathcal{H}$, if $2i-1 \in A_J$ then $\delta_i \leq Cn^{-(1\wedge\alpha_1)}$ for some constant $C > 0$; if $2i-1 \notin A_J$, $\delta_i \leq Cn^{-(1\wedge\alpha)}$ for some constant $C > 0$. This means

$$\frac{1}{n}\sum_i \delta_i^4 = \frac{1}{n}\sum_{i\in A_J}\delta_i^4 + \frac{1}{n}\sum_{i\notin A_J}\delta_i^4$$

$$\leq \frac{1}{n}M_{f3}n^{\gamma_f}Cn^{-4(1\wedge\alpha_1)} + \frac{1}{n}(n - M_{f3}n^{\gamma_f})Cn^{-4(1\wedge\alpha)}$$

$$= C_1 n^{-4(1\wedge\alpha)} + C_2 n^{\gamma_f - 1 - 4(1\wedge\alpha_1)} = O(n^{-(1\wedge 4\alpha)}).$$

On the other hand,

$$2\sum_{j,k}\gamma_{j,k}^2 = 2\sum_{j,k}\left(\sum_i W_{(j,k),i}(V_i - V(2i-1))\right.$$
$$\left. + \sum_i W_{(j,k),i}(V(2i-1) - \theta_{j,k})^2\right)$$

$$\leq 4\sum_{j,k}\left(\sum_i W_{(j,k),i}(V_i - V(2i-1))\right)^2$$

$$+ 4\sum_{j,k}\left(\sum_i W_{(j,k),i}V(2i-1) - \theta_{j,k}\right)^2$$

$$= 4\sum_i(V_i - V(2i-1))^2 + 4\sum_{j,k}\left(\sum_i W_{(j,k),i}V(2i-1) - \theta_{j,k}\right)^2.$$

Similarly to the previous calculation, we have

$$\sum_i(V_i - V(2i-1))^2 = O(n^{\gamma_V - 1 - 2(1\wedge\beta_1)} + n^{-2(1\wedge\beta)}).$$

It follows from Lemma 2 that

$$\sum_{j,k}\left(\sum_i W_{(j,k),i}V(2i-1) - \theta_{j,k}\right)^2$$

$$= \sum_k(\xi_{J/2,k} - V(2i-1))^2$$

$$= \sum_{k\in A_{J/2}}(\xi_{J/2,k} - V(2i-1))^2 + \sum_{k\notin A_{J/2}}(\xi_{J/2,k} - V(2i-1))^2$$

$$\leq C_1 n^{-2\beta} + C_2 n^{\gamma_V - 1 - 2\beta_1}.$$

The lemma is proved by putting these together. □

LEMMA 5. *Using the notation in Section 2.2, for any $x_* \in (0,1)$,*

$$\sup_{f\in\Lambda^\alpha(M_f,x_*,\delta), V\in\Lambda^\beta(M_V,x_*,\delta)}\left(\sum_k(\widetilde{\tau}_{j_0,k} - \xi_{j_0,k})\phi_{j_0,k}(x_*)\right.$$

$$\left. + \sum_{j,k}(\tau_{j,k} - \theta_{j,k})\psi_{j,k}(x_*)\right)^2$$

$$= O(n^{-(4\alpha\wedge 2\beta\wedge 1)}).$$

PROOF. It follows from the property of DWT that,

$$\left(\sum_k(\widetilde{\tau}_{j_0,k} - \xi_{j_0,k})\phi_{j_0,k}(x_*) + \sum_{j,k}(\tau_{j,k} - \theta_{j,k})\psi_{j,k}(x_*)\right)^2$$

$$= \left(\sum_i\left(\sqrt{\frac{2}{n}}\left(\frac{1}{2}\delta_i^2 + V_i\right) - \xi_{J-1,i}\right)\phi_{J-1,i}(x_*)\right)^2.$$

Note that $\phi(x)$ has compact support, say $\text{supp}(\phi) \subset [-L,L]$. So $\phi_{J-1,i}(x_* \neq 0)$ only if $\frac{2i}{n} \notin (x_* - \frac{2L}{n}, x_* + \frac{2L}{n})$. This means in the previous summation we only need to consider those $i$'s for which $\frac{2i}{n} \in (x_* - \frac{2L}{n}, x_* + \frac{2L}{n})$. For those $i$, $\text{supp}(\phi_{J-1,i}) \subset (x_* - \delta, x_* + \delta)$ for all sufficiently large $n$. On the interval $(x_* - \delta, x_* + \delta)$, both $f(x)$ and $V(x)$ has Lipschitz property and the lemma now follows from (31). □

Lemma 6 below shows that the estimator $\widehat{\text{Var}(D_i^2)}$ given in (20) has the desired property of being slightly positively biased.

LEMMA 6. *Suppose $V(x)$ is bounded away from zero and $z_i$'s are i.i.d. random variables. Suppose $\hat{\sigma}_k^2$ is the estimator mentioned in Section 2.3. Then for any $m > 0$ there exist constants $C_m > 0$ such that*

$$P\left(\bigcap_k \bigcap_{i \in \text{block } k} (E(e_i^2) < \hat{\sigma}_k^2 < 4E(e_i^2))\right) > 1 - C_m n^{-m}.$$

PROOF. Let $u_k$ denote the $k$th moment of $z_i$. It is easy to see that

$$E\Delta_i = V'_{2i-1} - V'_{2i} \leq 2M_V \left(\frac{1}{n}\right)^{\beta \wedge 1} + M_f \left(\frac{1}{n}\right)^{2(\alpha \wedge 1)},$$

$$E\Delta_i^2 = (V'_{2i-1} - V'_{2i})^2 + E(e_{2i-1}^2) + E(e_{2i}^2).$$

Since $E(e_i^2) = V_i^2(u_4 - 1) + 2\delta_i^2 V_i + 2\sqrt{2} V_i^{3/2} \delta_i u_3$, we know that

$$E(e_i^2) - E(e_j^2) \leq C_0\left(\left(\frac{|i-j|}{n}\right)^{\beta \wedge 1} + \left(\frac{|i-j|}{n}\right)^{\alpha \wedge 1}\right)$$

for some constant $C_0$. Denote by $B_k$ the set of indices in block $k$. Let $\omega_k = \max_{i \in B_k} \{E(e_i^2)\}$. Then for any $j \in B_k$

$$\omega_k - E(e_j^2) \leq C_0(n^{-(1-r)(\beta \wedge 1)} + n^{-(1-r)(\alpha \wedge 1)}) \leq C_0 n^{-(1-r)(\alpha \wedge \beta \wedge 1)}$$

and hence

$$E(\hat{\sigma}_k^2) = \frac{2}{(2 - 1/\log n)(n/2)^r} \sum_{2i \in B_k} ((V'_{2i-1} - V'_{2i})^2 + E(e_{2i-1}^2) + E(e_{2i}^2))$$

$$\geq \frac{2}{(2 - 1/\log n)(n/2)^r} \sum_{2i \in B_k} (E(e_{2i-1}^2) + E(e_{2i}^2))$$

$$\geq \frac{2}{(2 - 1/\log n)(n/2)^r} \sum_{2i \in B_k} (2\omega_k - 2C_0 n^{-(1-r)(\alpha \wedge \beta \wedge 1)})$$

$$= \omega_k + \frac{1/\log n}{2 - 1/\log n} \omega_k - \frac{2}{2 - 1/\log n} C_0 n^{-(1-r)(\alpha \wedge \beta \wedge 1)}.$$

Since $V(x)$ is bounded away from zero, we know that $\omega_k \geq C$ for some constant $C > 0$. So when $n$ is sufficiently large, there exists some constant $C_1 > 0$ such that
$$\frac{1/\log n}{2 - 1/\log n}\omega_k - \frac{2}{2 - 1/\log n}C_0 n^{-(1-r)(\alpha \wedge \beta \wedge 1)} > C_1/\log n.$$

Since all the moments of $e_i$ exist, all the moments of $\Delta_i$ exist. Then for any fixed positive integer $l$,

$$\begin{aligned}
P(\hat{\sigma}_k^2 &> \omega_k) \\
&= P(\hat{\sigma}_k^2 - E(\hat{\sigma}_k^2) > \omega_k - E(\hat{\sigma}_k^2)) \\
&\geq P\left(\hat{\sigma}_k^2 - E(\hat{\sigma}_k^2) > -\left(\frac{1/\log n}{2 - 1/\log n}\omega_k - \frac{2}{2 - 1/\log n}C_0 n^{-(1-r)(\alpha \wedge \beta \wedge 1)}\right)\right) \\
&\geq P(|\hat{\sigma}_k^2 - E(\hat{\sigma}_k^2)| < C_1/\log n) \\
&\geq 1 - \frac{E[(\hat{\sigma}_k^2 - E(\hat{\sigma}_k^2))^{2l}]}{(C_1/\log n)^{2l}} \\
&= 1 - \frac{1}{(C_1/\log n)^{2l} n^{2rl}(2 - 1/\log n)^{2l}} E\left(\sum_{2t \in B_k} (\Delta_t^2 - E\Delta_t^2)\right)^{2l}.
\end{aligned}$$

Since $\Delta_t$'s are independent random variables, we know that $E(\sum_{2t \in B_k}(\Delta_t^2 - E\Delta_t^2))^{2l}$ is of order $(n^r)^l$ for sufficiently large $n$. This means

$$(33) \qquad P(\hat{\sigma}_k^2 \geq \omega_k) = 1 - O\left(\frac{(\log n)^{2l}}{n^{rl}}\right).$$

So

$$P\left(\bigcap_{k=1}^{n^{1-r}} (\hat{\sigma}_k^2 > \omega_k)\right) \geq \left(1 - O\left(\left(\frac{\log^2 n}{n^r}\right)^l\right)\right)^{n^{1-r}} = 1 - O\left(\frac{(\log n)^{2l}}{n^{(l+1)r-1}}\right).$$

Since $l$ is an arbitrary positive integer, this means for any $m > 0$ there exists a constant $C_m > 0$ such that

$$P\left(\bigcap_k \bigcap_{i \in \text{block } k} (\hat{\sigma}_k^2 > E(e_i^2))\right) > 1 - C_m n^{-m}.$$

Similarly, we know that for any $m > 0$ there exists a constant $C_m' > 0$ such that

$$P\left(\bigcap_k \bigcap_{i \in \text{block}' k} (\hat{\sigma}_k^2 < 4E(e_i^2))\right) > 1 - C_m' n^{-m}. \qquad \square$$

A direct consequence of Lemma 6 is that $P(\bigcap_{j,k} \sigma_{j,k}^2 \leq \hat{\sigma}_{j,k}^2 \leq C\sigma_{j,k}^2) \geq 1 - C_m n^{-m}$ for any $m > 0$ and some constant $C_m$.

6.2. *Upper bound: Proof of Theorem 1.* It is clear that the estimators $\widehat{V}_e$, $\widehat{V}_0$, and thus $\widehat{V}$ have the same rate of convergence. Here we will only prove the convergence rate result for $\widehat{V}_e$. We shall write $\widehat{V}$ for $\widehat{V}_e$ in the proof. Note that

$$E\|\widehat{V} - V\|_{L_2}^2 = E\sum_{i=1}^{2^{j_0}}(\hat{\xi}_{j_0,i} - \xi_{j_0,i})^2 \qquad (34)$$
$$+ E\sum_{j=j_0}^{J_1}\sum_k(\widehat{\theta}_{j,k} - \theta_{j,k})^2 + \sum_{j=J_1+1}^{\infty}\sum_k \theta_{j,k}^2.$$

There are a fixed number of terms in the first sum on the RHS of (34). Equation (13) and Lemma 4 show that the empirical coefficients $\tilde{d}_{j_0,k}$ have variance of order $n^{-1}$ and sum of squared biases of order $O(n^{-(1\wedge 4\alpha \wedge 2\beta \wedge (1+2\beta_1 - \gamma_V))})$. Note that $\gamma_V - 1 - 2\beta_1 < \frac{2\beta}{1+2\beta}$, so

$$\sup_{V\in\Lambda^\beta(M)} E\sum_{i=1}^{2^{j_0}}(\hat{\xi}_{j_0,i} - \xi_{j_0,i})^2 = O(n^{-(1\wedge 4\alpha \wedge 2\beta \wedge (1+2\beta_1 - \gamma_V))}) + O(n^{-1})$$
$$= \max\left(O(n^{-4\alpha}), O\left(\frac{n}{\log n}\right)^{-2\beta/(1+2\beta)}\right).$$

Also, it is easy to see that the third sum on the RHS of (34) is small. Note that for $\theta_{j,k} = \langle V, \psi_{j,k}\rangle$, from Lemma 2,

$$\sum_{j=J_1+1}^{\infty}\sum_k \theta_{j,k}^2 = \sum_{j=J_1+1}^{\infty}\left(\sum_{k\in A_j}\theta_{j,k}^2 + \sum_{k\notin A_j}\theta_{j,k}^2\right)$$
$$\leq \sum_{j=J_1+1}^{\infty}(C_1 2^{j(\gamma_V - 1 - 2\beta_1)} + C_2 2^{-2j\beta})$$
$$= O\left(\left(\frac{n}{\log n}\right)^{-2\beta/(1+2\beta)}\right).$$

We now turn to the main term $E\sum_{j=j_0}^{J_1}\sum_k (\widehat{\theta}_{j,k} - \theta_{j,k})^2$. Note that

$$E\sum_{j=j_0}^{J_1}\sum_k(\widehat{\theta}_{j,k} - \theta_{j,k})^2 \leq 2E\sum_{j=j_0}^{J_1}\sum_k(\widehat{\theta}_{j,k} - \tau_{j,k})^2 + 2E\sum_{j=j_0}^{J_1}\sum_k(\tau_{j,k} - \theta_{j,k})^2.$$

The second term is controlled by Lemma 4. We now focus on the first term. Note that the thresholds $\lambda_{j,k}$ are random. We shall denote by $E_{|\lambda}(\cdot)$ the conditional expectation given all the thresholds $\lambda_{j,k}$. It follows from

Lemma 1 that
$$E(\hat{\theta}_{j,k} - \tau_{j,k})^2 = E(E_{|\lambda}(\hat{\theta}_{j,k} - \tau_{j,k})^2)$$
$$\leq E(\tau_{j,k}^2 \wedge 4\lambda_{j,k}^2) + E(z_{j,k}^2 I(|z_{j,k}| > \lambda_{j,k})).$$

Note that
$$E(z_{j,k}^2 I(|z_{j,k}| > \lambda_{j,k}))$$
$$= E(z_{j,k}^2 I(|z_{j,k}| > \lambda_{j,k}) I(\hat{\sigma}_{j,k}^2 \geq \sigma_{j,k}^2))$$
$$+ E(z_{j,k}^2 I(|z_{j,k}| > \lambda_{j,k}) I(\hat{\sigma}_{j,k}^2 \leq \sigma_{j,k}^2))$$
$$\leq E(z_{j,k}^2 I(|z_{j,k}| \geq \sigma_{j,k}\sqrt{2\log n})) + E(z_{j,k}^2 I(\hat{\sigma}_{j,k}^2 \leq \sigma_{j,k}^2))$$
$$\triangleq S_1 + S_2.$$

Set $\rho_{j,k} = \sigma_{j,k}\sqrt{2\log n}$. Since the moment generating functions of all $z_i^2$ exist in a neighborhood of the origin, there exists a constant $a > 0$ such that $E(e^{az_{j,k}/\sigma_{j,k}}) < \infty$. Let $A = C\log n$ for some constant $C > \max(1/a, 1)$, then
$$S_1 = E(z_{j,k}^2 I((A\sigma_{j,k}^2 \vee \rho_{j,k}^2) \geq |z_{j,k}| \geq \rho_{j,k}))$$
$$+ \sigma_{j,k}^2 E\left(\frac{z_{j,k}^2}{\sigma_{j,k}^2} I(|z_{j,k}| > (A\sigma_{j,k}^2 \vee \rho_{j,k}^2))\right)$$
$$\leq (A^2 \sigma_{j,k}^4 \vee \rho_{j,k}^4) P(|z_{j,k}| \geq \sigma_{j,k}\sqrt{2\log n})$$
$$+ \sigma_{j,k}^2 \frac{A^2 \vee 4(\log n)^2}{e^{a(A \vee 2\log n)}} E(e^{az_{j,k}/\sigma_{j,k}}).$$

Note that, when $2^j < n/(\log(n/2))^2$, each wavelet coefficient at level $j$ is a linear combination of $m \geq (\log n)^2$ of the $y_i$'s. It then follows from Lemma 3 and (32) that
$$P(|z_{j,k}| > \sigma_{j,k}\sqrt{2\log(n/2)}) \leq O(n^{-(1-a)})$$
for any $a > 0$. Also,
$$\frac{A^2 \vee 4(\log n)^2}{e^{a(A \vee 2\log n)}} \leq \frac{C^2 \log^2 n}{e^{aC\log n}} \leq O\left(\frac{\log^2 n}{n}\right).$$

Combining these together, and since $\sigma_{j,k}^2 = O(1/n)$, we have
$$S_1 \leq O\left(\frac{\log^2 n}{n^2}\right).$$

This means $S_1$ is negligible as compared to the upper bound given in (25).

It is easy to see that $S_2 \leq (E(z_{j,k}^4) P(\hat{\sigma}_{j,k}^2 \leq \sigma_{j,k}^2))^{1/2}$. Lemma 6 yields $P(\hat{\sigma}_{j,k}^2 \leq \sigma_{j,k}^2) = O(n^{-m})$ for any $m > 0$. So $S_2$ is also negligible.

We now turn to $E(\tau_{j,k}^2 \wedge 4\lambda_{j,k}^2)$. Note that

(35) $$E(\tau_{j,k}^2 \wedge 4\lambda_{j,k}^2) \le 2(\tau_{j,k} - \theta_{j,k})^2 + 2E(\theta_{j,k}^2 \wedge 4\lambda_{j,k}^2).$$

The first part is controlled by Lemma 4. For the second part,

$$E(\theta_{j,k}^2 \wedge 4\lambda_{j,k}^2) \le E(\theta_{j,k}^2 \wedge (4 \times 4\rho_{j,k}^2)) + E(\theta_{j,k}^2 I(\hat\sigma_{j,k}^2 > 4\sigma_{j,k}^2))$$
$$\le 4(\theta_{j,k}^2 \wedge 4\rho_{j,k}^2) + \theta_{j,k}^2 P(\hat\sigma_{j,k}^2 > 4\sigma_{j,k}^2).$$

Note that $\theta_{j,k}^2$ is bounded (see Lemma 2). From Lemma 6, $P(\hat\sigma_{j,k}^2 > 4\sigma_{j,k}^2) = O(n^{-m})$ for any $m > 0$. So $\theta_{j,k}^2 P(\hat\sigma_{j,k}^2 > 4\sigma_{j,k}^2)$ is negligible as compared to the upper bound in (25).

We now turn to $\theta_{j,k}^2 \wedge 4\rho_{j,k}^2$, note that $j \ge \frac{J-\log_2 J}{2\beta+1}$ implies $2^j \ge (\frac{n}{\log n})^{1/(2\beta+1)}$, and

$$\begin{cases} \theta_{j,k}^2 \wedge 4\rho_{j,k}^2 \le 4\rho_{j,k}^2 \le C\left(\dfrac{\log n}{n}\right), & \text{if } j \le \dfrac{J-\log_2 J}{2\beta+1} \text{ and } k \notin A_j, \\ \theta_{j,k}^2 \wedge 4\rho_{j,k}^2 \le \theta_{j,k}^2 \le C 2^{-j(1+2\beta)}, & \text{if } j \ge \dfrac{J-\log_2 J}{2\beta+1} \text{ and } k \notin A_j, \\ \theta_{j,k}^2 \wedge 4\rho_{j,k}^2 \le \theta_{j,k}^2 \le C 2^{-j(1+2\beta_1)}, & \text{if } k \in A_j. \end{cases}$$

This means

$$\sum_{j,k} \theta_{j,k}^2 \wedge 4\rho_{j,k}^2 \le \sum_{j \le (J-\log_2 J)/(2\beta+1)} C 2^j \left(\frac{\log n}{n}\right) + \sum_{j > (J-\log_2 J)/(2\beta+1)} C 2^{-j 2\beta}$$
$$+ \sum_j \sum_{k \in A_j} C 2^{-j(1+2\beta_1)}$$
$$\le C\left(\frac{\log n}{n}\right) 2^{(J-\log_2 J)/(2\beta+1)}$$
$$+ C 2^{-2\beta(J-\log_2 J)/(2\beta+1)} + C 2^{\gamma_V - 1 - 2\beta_1}$$
$$\le C\left(\frac{\log n}{n}\right)^{2\beta/(1+2\beta)}.$$

Putting the above bounds together, one can easily see that

$$\sum_{(j,k)} E((\tau_{j,k})^2 \wedge 4\lambda_{j,k}^2) \le M_f^4 n^{-4\alpha} + 4M_V^2 (n^{-2} \wedge n^{-2\beta}) + C\left(\frac{n}{\log n}\right)^{-2\beta/(1+2\beta)}$$
$$\le C \max\left(n^{-4\alpha}, \left(\frac{n}{\log n}\right)^{-2\beta/(1+2\beta)}\right).$$

This proves the global upper bound (25).

6.3. *Upper bound: Proof of Theorem 2.* We now consider the bound given in (26) under pointwise MSE. Without loss of generality, we shall assume that $f$ and $V$ are in the Lipschitz classes instead of the local Lipschitz classes, that is, we assume $f \in \Lambda^\alpha(M_f)$ and $V \in \Lambda^\beta(M_V)$. Note that

$$E(\widehat{V}(x_*) - V(x_*))^2$$
$$= E\left(\sum_{i=1}^{2^{j_0}}(\hat{\xi}_{j_0,i} - \xi_{j_0,i})^2 \phi_{j_0,k}(x_*)\right.$$
$$\left. + \sum_{j=j_0}^{J_1}\sum_k (\hat{\theta}_{j,k} - \theta_{j,k})\psi_{j,k}(x_*) + \sum_{j>J_1,k} \theta_{j,k}\psi_{j,k}(x_*)\right)^2$$
$$\leq 3\left(\sum_{i=1}^{2^{j_0}}(\hat{\xi}_{j_0,i} - \xi_{j_0,i})^2 \phi_{j_0,k}(x_*)\right)^2 + 3\left(\sum_{j=j_0}^{J_1}\sum_k (\hat{\theta}_{j,k} - \theta_{j,k})\psi_{j,k}(x_*)\right)^2$$
$$+ 3\left(\sum_{j>J_1}\sum_k \theta_{j,k}\psi_{j,k}(x_*)\right)^2$$
$$\triangleq I_1 + I_2 + I_3.$$

$I_1$ is bounded in the same way as in the global case. Since we are using wavelets of compact support, there are at most $L$ basis functions $\psi_{j,k}$ at each resolution level $j$ that are nonvanishing at $x_*$ where $L$ is the length of the support of the wavelet $\psi$. Denote $K(j, x_*) = \{k : \psi_{j,k}(x_*) \neq 0\}$. Then $|K(j, x_*)| \leq L$. Hence

$$I_3 = 3\left(\sum_{j>J_1}\sum_{k \in K(j,x_*)} \theta_{j,k}\psi_{j,k}(x_*)\right)^2 \leq 3\left(\sum_{j>J_1} CL2^{-j(\beta+1/2)}2^{j/2}\right)^2$$
$$= O(2^{-J_1\beta}) = o(n^{-2\beta/(1+2\beta)}).$$

We now turn to $I_2$. First,

$$I_2 \leq 3\left(\sum_{j,k}(E(\hat{\theta}_{j,k} - \theta_{j,k})^2)^{1/2}|\psi_{j,k}(x_*)|\right)^2.$$

Note that

$$E(\hat{\theta}_{j,k} - \theta_{j,k})^2 \leq 2E(\hat{\theta}_{j,k} - \tau_{j,k})^2 + 2(\tau_{j,k} - \theta_{j,k})^2$$
$$\leq 4(\tau_{j,k} - \theta_{j,k})^2 + 2E(\theta_{j,k}^2 \wedge 4\lambda_{j,k}^2) + 2E(z_{j,k}^2 I(|z_{j,k}| > \lambda_{j,k}))$$
$$\leq 4(\tau_{j,k} - \theta_{j,k})^2 + 8(\theta_{j,k}^2 \wedge 4\rho_{j,k}^2) + 2\theta_{j,k}^2 P(\hat{\sigma}_{j,k}^2 > 4\sigma_{j,k}^2)$$
$$+ 2E(z_{j,k}^2 I(|z_{j,k}| > \lambda_{j,k})).$$

This means

$$I_2 \leq 96\left(\sum_{j,k}(\theta_{j,k}^2 \wedge 4\rho_{j,k}^2)^{1/2}\psi_{j,k}(x_*)\right)^2$$
$$+ 48\left(\sum_{j,k}(\tau_{j,k} - \theta_{j,k})\psi_{j,k}(x_*)\right)^2$$
$$+ 24((\theta_{j,k}^2 P(\hat{\sigma}_{j,k}^2 > 4\sigma_{j,k}^2))^{1/2}\psi_{j,k}(x_*))^2$$
$$+ 24((E(z_{j,k}^2 I(|z_{j,k}| > \lambda_{j,k}))\psi_{j,k}(x_*))^{1/2})^2.$$

The last two terms follow from the proof of the global upper bound and the second term is controlled by Lemma 5. For the first term, from the discussion in the proof of global upper bound, we have

$$\sum_{j,k}(\theta_{j,k}^2 \wedge 4\rho_{j,k}^2)^{1/2}\psi_{j,k}(x_*)$$
$$= \sum_{j_0 \leq j \leq (J-\log_2 J)/(2\beta+1)} \sum_{k \in K(j,x_*)} (\theta_{j,k}^2 \wedge 4\rho_{j,k}^2)^{1/2}\psi_{j,k}(x_*)$$
$$+ \sum_{j > (J-\log_2 J)/(2\beta+1)} \sum_{k \in K(j,x_*)} (\theta_{j,k}^2 \wedge 4\rho_{j,k}^2)^{1/2}\psi_{j,k}(x_*)$$
$$\leq \sum_{j_0 \leq j \leq (J-\log_2 J)/(2\beta+1)} CL2^{j/2}\left(\frac{\log n}{n}\right)^{1/2}$$
$$+ \sum_{j > (J-\log_2 J)/(2\beta+1)} CL2^{j/2}2^{-j(\beta+1/2)}$$
$$= O\left(\left(\frac{\log n}{n}\right)^{\beta/(1+2\beta)}\right).$$

Putting these together, one can see that $I_2 \leq C\max(n^{-4\alpha}, (\frac{n}{\log n})^{2\beta/(1+2\beta)})$. This proves the local upper bound (26).

6.4. *Lower bound: Proof of Theorem 3.* We first outline the main ideas. The constrained risk inequality of Brown and Low (1996) implies that if an estimator has a small risk $\varepsilon^2$ at one parameter value $\theta_0$ and $(\theta_1 - \theta_0)^2 \gg \varepsilon\rho$ where $\rho$ is the chi-square affinity between the distributions of the data under $\theta_0$ and $\theta_1$, then its risk at $\theta_1$ must be "large." Now the assumption (27) means that the estimator $\widehat{V}(x_*)$ has a small risk at $\theta_0 = V_0(x_*)$. If we can construct a sequence of functions $V_n$ such that $V_n$ is "close" to $V_0$ in the sense that $\rho$ is small and at the same time $\Delta = |V_n(x_*) - V_0(x_*)|$ is "large," then it follows from the constrained risk inequality that $\widehat{V}(x_*)$ must have a "large" risk at $\theta_1 = V_n(x_*)$. So the first step of the proof is a construction for such a sequence of functions $V_n$.

Set $V_0 \equiv 1$ and let $g$ be a compactly supported, infinitely differentiable function such that $g(0) > 0$, $\int g = 0$ and $\int g^2 = 1$. Set
$$V_n(x) = V_0(x) + \tau_n^\beta g(\tau_n^{-1}(x - x_*)),$$
where $\tau_n = (\frac{c \log n}{n})^{1/(1+2\beta)}$ and $0 < c \leq 1$ is a constant. It is easy to check that $f_n$ are in $\Lambda^\beta(M)$ if the constant $c$ is chosen sufficiently small.

The chi-square affinity can be bounded same as before. Note that the chi-square affinity between $\Phi = N(0,1)$ and $\Psi = N(0, 1+\gamma_n)$ is

(36) $$\rho(\Phi, \Psi) = (1 - \gamma_n^2)^{-1/2}.$$

Let $y_i = V(x_i)z_i$, $i = 1, \ldots, n$, where $z_i$ are i.i.d. $N(0,1)$ variables. Denote by $P_0$ and $P_n$ the joint distributions of $y_1, \ldots, y_n$ under $V = V_0$ and $V = V_n$, respectively. Then it follows from (36) that

$$\begin{aligned}
\rho_n &\equiv \rho(P_0, P_n) \\
&= \prod_{i=1}^n [1 - \tau_n^{2\beta} g^2(\tau_n^{-1}(x_i - x_*))]^{-1/2} \\
&= \exp\left\{-\frac{1}{2} \sum_{i=1}^n \log(1 - \tau_n^{2\beta} g^2(\tau_n^{-1}(x_i - x_*)))\right\} \\
&\leq \exp\left\{\tau_n^{2\beta} \sum_{i=1}^n g^2(\tau_n^{-1}(x_i - x_*))\right\}
\end{aligned}$$

where the last step follows from the fact $-\frac{1}{2}\log(1-z) \leq z$ for $0 < z < \frac{1}{2}$. Note that $(n\tau_n)^{-1} \sum_{i=1}^n g^2(\tau_n^{-1}(x_i - x_*)) \to \int g^2 = 1$, so $\sum_{i=1}^n g^2(\tau_n^{-1}(x_i - x_*)) \leq 2n\tau_n$ for sufficiently large $n$, and hence $\rho_n \leq \exp(2n\tau_n^{1+2\beta}) \leq n^{2c}$. Since the zero function is in $\Lambda^{\alpha_0}(M_f)$ and $V_0 \in \Lambda^{\beta_0}(M_V)$, equation (27) implies that for some constant $C > 0$,
$$E(\widehat{V}(x_*) - V_0(x_*))^2 \leq Cn^{-2\beta_1/(1+2\beta_1)} n^{-r}$$
for $r = \min\{4\alpha_0, \frac{2\beta_0}{1+2\beta_0}\} - \frac{2\beta_1}{1+2\beta_1} > 0$. Hence, for sufficiently large $n$, it follows from the constrained risk inequality of Brown and Low (1996) that

$$\begin{aligned}
E(\widehat{V}(x_*) - V_n(x_*))^2 &\geq \tau_n^{2\beta} g^2(0)\left(1 - \frac{2C^{1/2} n^{-\beta_1/(1+2\beta_1)} n^{-r/2} n^c}{\tau_n^\beta g(0)}\right) \\
&= \left(\frac{c \log n}{n}\right)^{2\beta/(1+2\beta)} \\
&\quad \times g^2(0)\left(1 - \frac{2C^{1/2} n^{-\beta_1/(1+2\beta_1)} n^{-r/2} n^c}{(c\log n/n)^{\beta_1/(1+2\beta_1)} g(0)}\right)
\end{aligned}$$

$$\geq \frac{1}{2} c^{2\beta_1/(1+2\beta_1)} g^2(0) \cdot \left(\frac{\log n}{n}\right)^{2\beta_1/(1+2\beta_1)}$$

by choosing the constant $c \leq r/2$.

**Acknowledgments.** We wish to thank the Associate Editor and two referees for their constructive comments which led to an improvement in some of our earlier results and also helped with the presentation of the paper.

Department of Statistics
The Wharton School
University of Pennsylvania
Philadelphia, Pennsylvania 19104
USA
E-mail: tcai@wharton.upenn.edu
liewang@wharton.upenn.edu